\documentclass{article}

\usepackage{latexsym}
\usepackage{amssymb}
\usepackage{amscd,amssymb,amsmath}
\usepackage{graphicx}
\usepackage{color}
\usepackage{setspace}
\usepackage{float}
\usepackage{indentfirst}
\usepackage{tabularx}
\usepackage{longtable}
\usepackage[all]{xy}

\newtheorem{theorem}{Theorem}[section]
\newtheorem{definition}[theorem]{Definition}
\newtheorem{proposition}[theorem]{Proposition}

\newtheorem{lemma}[theorem]{Lemma}
\newtheorem{remark}[theorem]{Remark}
\newtheorem{example}[theorem]{Example}
\def\0{\underline 0}

\newenvironment{prova}{\noindent{\it Proof: }}{\hfill $\square$}
\begin{document}

%
%

\title{Results on Milnor Fibrations for mixed polynomials with non-isolated singularities}

\author{N. G. Grulha Jr.\footnote{Instituto de Ci\^encias Matem\'aticas e de Computa\c{c}\~ao - USP. Av. Trabalhador s\~ao-carlense, 400 - Centro.
  CEP: 13566-590 - S\~ao Carlos - SP, Brazil.
  njunior@icmc.usp.br} \and 
R. S. Martins \footnote{Universidade Federal de Santa Catarina - UFSC.
 R. Dona Francisca, 8300 - Distrito Industrial.
CEP: 89219-600 - Joinville - SC, Brazil.
rafaella.sm@ufsc.br}
}

\maketitle


\setcounter{section}{0}

\begin{abstract}
\noindent In this article we investigate mixed polynomials and present conditions that can be applied on a specific class of polynomials in order to prove the existence of the Milnor Fibration, Milnor-Lê Fibration and the equivalence between them. We prove for this class the of functions that the Milnor-L\^e fiber on a regular value is homeomorphic to the Milnor-L\^e fiber on a critical value. We develop a criterion to verify the transversality property and apply it to a special case of the class of mixed polynomial.

\vspace{0.8cm}


\vspace{0.4cm}

\end{abstract}

\section*{Introduction}
\label{sec:1}
In 1960's Milnor proved the result, nowadays known as Milnor's Fibration Theorem (see Milnor \cite{Milnor}). The theorem concerns the geometry and topology of holomorphic map-germs $f:(\mathbb{C}^{n},0)\to (\mathbb{C},0)$. This result is a landmark in singularity theory with a huge impact in many branches in mathematics.

When the map-germ in question has an isolated critical point, the fiber of the Milnor Fibration is diffeomorphic to a $2(n-1)$-ball to which one attaches handles of middle index. The number of such handles is what we call the Milnor number of the singularity. If the critical point is non-isolated, then the Milnor fiber is diffeomorphic to a ball to which we must attach handles of various indices, in this direction we have, for instance the works of Massey and Tib\u{a}r \cite{Massey1,Massey2,Tibar}.

Milnor proved that in the isolated singular case these two fibrations described below are equivalent.

If we set $V= f^{-1}(0)$ and  $K_\varepsilon = \mathbb{S}_\varepsilon \cap   V$ where $\mathbb{S}_\varepsilon$ is the sphere of radius $\varepsilon$ in $\mathbb{C}^n$ centered at $0$, then the first Milnor Fibration is
$$ \phi: \frac{f}{|f|}: \mathbb{S}_\varepsilon  \setminus K_\varepsilon \longrightarrow \mathbb{S}^1 \ , $$
for $\varepsilon  > 0$ sufficiently small (see Milnor, \cite{Milnor}). Given $\varepsilon$ as above,
the second fibration is
$$ \frac{f}{|f|}: N(\varepsilon,\delta) \longrightarrow \mathbb{S}^1 \,,$$
where $N(\varepsilon,\delta)$ is the tube $f^{-1}(\partial\mathbb{D}_\delta) \cap \mathbb{B}_\varepsilon$ for $\varepsilon \gg \delta > 0$, with  $\mathbb{B}_\varepsilon $ being the ball in $\mathbb{C}^n$ bounded by $\mathbb{S}_\varepsilon = \partial \mathbb{B}_\varepsilon $ and $\partial\mathbb{D}_\delta$ the boundary of the disc $\mathbb{D}_\delta$ in $\mathbb{C}$ (see L\^e, \cite{Le}).

The first of these is known as the Milnor Fibration of the map-germ, while the fibration on the tube is called nowadays as Milnor-L\^e Fibration.

The literature on this topic is vast and we refer for instance to \cite{Seade-50} for an updated survey paper by Seade.

The equivalent theory for real analytic map-germs also springs from \cite{Milnor} but it is much less developed. After some important articles by various authors, as for instance by Looijenga \cite{Loo}, Church and Lamotke \cite{CL}, A’Campo \cite{ACampo}, Perron
\cite{Perron}, Kauffman and Neumann \cite{KN}, Jacquemard \cite{Jac} and others, in the 1970s and early 1980s the subject became somehow dormant till the 1990s, when Ruas, Seade and Verjovsky in \cite {Ruas-Seade-Verjovsky} opened new lines of research. For example, Ara\'ujo, Cisneros-Molina, Dutertre, Menegon, Snoussi and Seade in \cite{NA,CMSS,N} investigate the non-isolated singularities in the real case.

In a more general situation, the study of Milnor Fibrations defined by functions $f\overline{g}$ essentially began with the works of Pichon and Seade \cite{Pichon,Pichon-Seade2}. In a different approach Ara\'{u}jo, Ribeiro and Tib\u{a}r in \cite{ART} also study the Milnor \linebreak Fibration on real analytc maps.
In \cite{CGS} Cisneros-Molina, Grulha and Seade study the topology of the fibers of real analytic maps $f: \mathbb{R}^n \to \mathbb{R}^p , n > p$, in a neighborhood of a critical point.

In \cite{Ruas-Seade-Verjovsky}, followed by \cite{Seade-libro},  the authors consider   germs at the origin in $U \subset \mathbb{C}^n$ of holomorphic vector fields $X(z_1, \cdots, z_n)=(t_1 z_1^{b_1}, \cdots, t_n z_n^{b_n})$ and $F(z_1, \cdots, z_n)=(k_1 z_1^{a_1}, \cdots, k_n z_n^{a_n})$,  where $\sigma$ is some permutation of $\{1,...,n\}$, $k_j, t_j \in \mathbb{C}^*$, and  $a_j, b_j \in \mathbb{Z}^+$, and  their Hermitian product $\psi: U \subset \mathbb{C}^n \longrightarrow \mathbb{C}$:
\begin{equation}\label{psi-introducao}
\psi(z):= <F(z),X(z)> = \sum^{n}_{i=1}F_i(z) \overline{X_i}(z) = \sum_{j=1}^n k_j \overline{t_j} z_j^{a_j}\overline{z_j}^{b_j}.
\end{equation}

It turns out that if one has  $a_j \ne b_j$ for all $j$, then $\psi$ has a isolated critical point at $0$ and has Milnor Fibrations, as in the holomorphic case, on sphere and tube.

The proof of this result uses the fact that these singularities are a reminiscent of the classical weighted homogeneous ones, since $\psi^{-1}(0)$ admits a canonical action of the Lie group $\mathbb{S}^1 \times \mathbb{R}^+$. The equivalence between these fibrations (on the sphere and on the tube) is given by the flow of the $\mathbb{R}^+$-action. Furthermore, it is proved in \cite {Ruas-Seade-Verjovsky, Oka-dife} that the corresponding Milnor Fibrations for $\psi$ described in (\ref{psi-introducao}) actually are smoothly equivalent to that of the classical Pham-Brieskorn singularities.

 This gave rise to the theory of mixed and polar weighted singularities, developed by M. Oka, Cisneros-Molina and others (see \cite{Seade-50} for an account on the subject). Mixed singularities are defined by complex valued real analytic functions on the variables $(z_1,\cdots, z_n, \bar z_1, \cdots, \bar z_n)$ and M. Oka proved in \cite{Oka} important results about Milnor Fibrations for these singularities proving they are ``strongly non-degenerate", a condition that springs from looking at the associated Newton polygon, as in the holomorphic setting.

In this paper we consider the complementary case of map-germs.
\begin{equation}\label{psi-introducao1}
\psi(z):= <F(z),X(z)> = \sum^{n}_{i=1}F_i(z) \overline{X_i}(z) = \sum_{j=1}^n k_j \overline{t_j} z_j^{a_j}\overline{z_j}^{b_j},
\end{equation}
where  $a_j = b_j$ for at least one $j$, $j = 1, \cdots,n$. In this situation we investigate,
separately, the cases  $a_j = b_j$ for all $j$, $j = 1, \cdots,n$ and $a_j \neq b_j$ for at least one $j$, $j = 1, \cdots,n$.

These type of functions are not polar weighted because the condition of \linebreak having $a_j = b_j$ for at least one $j$ breaks the $S^1$-action, so the techniques from \linebreak \cite{Ruas-Seade-Verjovsky, Seade-libro}  does not apply either. Moreover, these singularities are not strongly non-degenerate, so Oka's theorem does not apply either. But our results are also given by the combinatorial information obtained from the mixed polynomials.

In Section 1, we present some general results to determine when a real \linebreak analytic map has the transversality property (Definition \ref{defproptransversalidade}), which is equivalent to admit a Milnor-L\^e Fibration on a tube. We also give examples which show interesting features regarding the behavior of the critical points that motivate our main results, Theorem \ref{teofibracaoitens} and Proposition \ref{teoremaespecial}.

In Section 2, we first show that the map $\psi$ that we envisage here are all weighted homogeneous with an $\mathbb R^+$-action. This implies that their discriminant is linear. In fact, Proposition  \ref{caracdodiscriminante} says that it is generically the union of line-segments. This also implies that once we have a Milnor-L\^e Fibration on a tube, the  $\mathbb R^+$-action carries this fibration onto a Milnor Fibration on the sphere, and therfore these singularities are $d$-regular (see Definition \ref{definicaod-regular}).

In Section 3, we show that there are two fibrations for a special case of $\psi$, using the transversality property criterion developed in section 1. Then we show that one has three essentially different possibilities for the type of critical points of $\psi$ depending on the exponents $a_j, b_j$ which are equal. A careful analysis of these cases yields to Theorem \ref{teofibracaoitens}.


\section{A criterium for local Milnor-L\^e Fibrations}\label{sec:2}

Let us consider the real analytic map-germ $f: (\mathbb R^n,\0) \longrightarrow (\mathbb R^p,0)$, where  $n > p \ge 1$, with a critical point at the origin $\0 \in \mathbb R^n$, such that $V:= f^{-1}(0)$ has dimension greater then zero.

Recall that a local Milnor sphere for $f$ means a sphere in $\mathbb R^n$ of radius $\varepsilon_0 > 0$ and center at $\0$, such that every sphere  in $\mathbb R^n$ of  radius less then $\varepsilon_0$ intersects $V$ transversally with respect to some Whitney stratification. Such spheres always exist by \cite{Milnor, BV}.

Let us recall that $f$ admits a local Milnor-L\^e Fibration \cite{Le} if for every sufficiently small Milnor sphere $\mathbb{S}_\varepsilon$, boundary of the  ball $\mathbb B_\varepsilon$, there exists
$\delta >0$ sufficiently small with respect to $\varepsilon$, where $\mathbb D_\delta$ denotes the ball in $\mathbb R^p$ of center at $0$ and radius $\delta$, the image of the critical points of $f$ denoted by  $\Delta_f \subset \mathbb R^p$,
and $N(\varepsilon,\delta) := \big(f^{-1}(\mathbb D_\delta \setminus \Delta_f)\big) \cap \mathbb B_\varepsilon$, we have

$$f: N(\varepsilon,\delta) \longrightarrow \mathbb D_\delta \setminus \Delta_f, $$
is a locally trivial fibration. We call $\Delta_f$ the discriminant of $f$.

In \cite{CMSS} one finds the following definition, which is a necessary and sufficient condition for having a local Milnor-L\^e Fibration:

\begin{definition}\label{defproptransversalidade} {\rm A map-germ $f$ as above has} the transversality property at $\0$ {\rm if there exists a real number $\varepsilon_0$ with $\varepsilon_0 > 0$ such that, for every $\varepsilon$ with $0 < \varepsilon \leq \varepsilon_0$, there exists a real number $\delta= \delta(\varepsilon)$, with $0 < \delta \ll \varepsilon$, such that for every $t \in \mathbb{B}_\delta^p \setminus \Delta_f$ one has that either $f^{-1}(t)$ does not intersect the sphere $\mathbb{S}_\varepsilon^{n-1}$ or $f^{-1}(t)$ intersects $\mathbb{S}_\varepsilon^{n-1}$ transversally in $\mathbb{R}^n$.}
\end{definition}

\begin{remark} Notice that in general the image of a map germ and its discri\-minant may not be well defined as set germs. For a descriptive analyzes on the situation see \cite{ART}.

\end{remark}

In this section we give a computable condition to verify if a map-germ $f$ has the transversality property. In order to do so,
we state next proposition that classify the points where the fibers of $f$ are transversal to the spheres, that is, the points  where $T_xf^{-1}(f(x)) \nsubseteq T_x \mathbb{S}_\varepsilon^{n-1}$ for $x \in f^{-1}(f(x)) \cap \mathbb{S}_\varepsilon^{n-1}$.

\begin{proposition}\label{propdastangencias}
Let $f:(\mathbb{R}^n,0) \longrightarrow (\mathbb{R}^p,0)$, with $n> p \geq 1$, be a map-germ of class $\mathcal{C}^\ell, \ \ell \geq 1$. Consider $x \in f^{-1}(f(x)) \cap \mathbb{S}^{n-1}_\varepsilon$ then $f^{-1}(f(x))$ and $\mathbb{S}^{n-1}_\varepsilon$ are transversal if and only if the vectors $x$ and the gradients of coordinate functions of $f$ at $x$ are linearly independent.
\end{proposition}

\begin{prova}
Consider $x \in f^{-1}(f(x))$, where $f(x) \notin \Delta_f$ and $x \in f^{-1}(f(x)) \cap \mathbb{S}_\epsilon^{n-1}$,\linebreak with $\epsilon > 0$ sufficiently small and let us suppose
\begin{equation}\label{xlinear}
 x = \alpha_1 \nabla f_1(x) + \cdots + \alpha_p \nabla f_p(x), \ \ \alpha_1, \cdots, \alpha_p \in \mathbb{R}.
\end{equation}
To have $f^{-1}(f(x))$ and $\mathbb{S}^{n-1}_\epsilon$ transversal one has that
$$T_x f^{-1}(f(x)) + T_x\mathbb{S}_\epsilon^{n-1} = \mathbb{R}^n,$$
where $T_x \mathbb{S}^{n-1}_\epsilon = \{x\}^\bot$ has dimension $n-p$ and $T_x f^{-1}(f(x)) = Ker \ df(x)$ of dimension $n-2$. Transversity does not occur when $T_x f^{-1}(f(x)) \subset T_x \mathbb{S}^{n-1}_\epsilon$.

Let $v \in T_x f^{-1}(f(x))$, in this case we have $<\nabla f_1(x), v> =0, \cdots, \linebreak <\nabla f_p(x), v> =0$, using (\ref{xlinear}) it follows
\begin{eqnarray*}
<\nabla f_1(x), v> + \cdots + <\nabla f_p(x),v> & = & 0  \\
\alpha_1 <\nabla f_1(x), v> + \cdots + \alpha_p <\nabla f_p(x),v> &  = & 0\\
<\alpha_1 \nabla f_1(x) + \cdots + \alpha_p \nabla f_p(x),v>& = & 0 \\
<x,v> & = & 0  \  \Leftrightarrow v \in T_x\mathbb{S}^{n-1}_\epsilon.
\end{eqnarray*}

Therefore, $T_x f^{-1}(f(x)) \subset T_x \mathbb{S}^{n-1}_\epsilon$ and we can conclude that at the point $x$ the manifold $f^{-1}(f(x))$ and $\mathbb{S}^{n-1}_\epsilon$ are not transversal.

\end{prova}

The next example is due to Snoussi. The corresponding map does not have the transversality property.

\begin{example}\label{exfnaotransversal}{\rm Consider the real analytic map
$$
\begin{array}{cccc}
f \ : & \! \mathbb{R}^3  & \! \longrightarrow
& \! \mathbb{R}^{2}\\
& \! (x,y,z) & \! \longmapsto
& \! (xy+z^2,x).
\end{array}
$$
The set of zeros of $f$ is given by $V = \{(0,y,0); \ y \in \mathbb{R}\}$. The jacobian matrix is
$$ Jf = \left[
\begin{array}{cccccc}
y & x & 2z \\
1 & 0 & 0 \\
\end{array}
\right],$$
the critical set is $\{(0,y,0); \ y \in \mathbb{R}\}$, hence the discriminant is $\Delta_f = \{(0,0)\}$ and $f$ has an isolated critical value.

Let us see whether  $f$ has the transversality property.
The tangent space of the fiber $f^{-1}(f(p))$ at $p$ translated to the origin is the line $\{(0,0,z) : z \in \mathbb{R} \}$, defined by the equations
$$
\left\{
\begin{array}{lcl}
bx+ay & = &  0\\
x & = &  0
\end{array}
\right.
$$

The tangent space at $p$ of the sphere in $\mathbb{R}^3$ centered on the origin containing $p$ translated to the origin is given by the equation
$$ax+by=0.$$

Hence the tangent line to the fiber of $f$ is contained in the tangent space of the sphere, consequently they are not transverse. Taking the parameterized curve $\gamma(a) = (a, \sqrt{\varepsilon^2-a^2},0)$ on $\mathbb{S}^2_\varepsilon$ then $f(\gamma(a)) = (a\sqrt{\varepsilon^2-a^2},a)$. If we take $a$ tends to zero we have the curve $\gamma(a)$ in the sphere $\mathbb{S}_\varepsilon^2$, approaching $V \cap \mathbb{S}_\varepsilon^2$ and such that $\displaystyle\lim_{a \to 0} f(\gamma(a)) = 0$.
This holds for all such spheres $\mathbb{S}_\varepsilon^2$, so $f$ does not have the transversality property.

Another way to verify that we do not have the transversality property  is to observe that $f$ has an isolated critical value, then if we had a locally trivial fibration away from $\{0\}$, the topology of all fibers outside $V$ would be the same. Yet, the fiber over $(1,0)$ is a disconnected manifold while the fiber over $(0,1)$ is connected.

In particular, the map $f$ does not have the Thom $a_f$ - Property (Defini\-tion 3.1, \cite{CSS-dr}), though one can check that it is $d$-regular by Proposition 3.8 of \cite{CSS-dr}.
}
\end{example}

Inspired by the above example, the next proposition gives us a criterion to have the transversality property.

\begin{proposition}\label{teopropdetransversalidade}
Consider a map-germ $f:(\mathbb{R}^n,0) \longrightarrow (\mathbb{R}^p,0)$, $n> p \geq 1$,  of class $\mathcal{C}^\ell, \ell \geq 1$, and consider a local Milnor sphere $\mathbb{S}_\varepsilon^{n-1} \subset \mathbb{R}^n$ for $f$. Then $f$ does not have the transversality property with respect to that sphere $\mathbb{S}_\varepsilon^{n-1}$
if and only if there is a sequence $(x_n)_{n \in \mathbb{N}}$ in $\mathbb{S}_\varepsilon^{n-1}$ satisfying
$\displaystyle\lim_{n\rightarrow\infty} x_n = p \in V \cap \mathbb{S}_\varepsilon^{n-1}$ and such that
 the vectors $x_n$ and the gradients of $f$ at  $x$ are linearly dependent.
  \end{proposition}

Notice that $f$ fails to have the transversality property if and only if there are spheres as above of arbitrarily small radius.

\begin{figure}[h]
\begin{center}
\includegraphics[scale=0.5,angle=0]{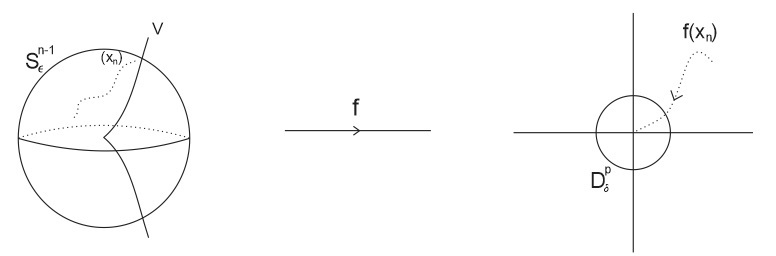}
 \caption{The sequence that loses the transversality property of the $f$}
  \end{center}
  \end{figure}

In Example \ref{ex15} the map g is given by the mixed polynomial $\psi(z)= \linebreak  z_1\overline{z_1} + z_2^2\overline{z_2}$ written in real coordinates and $\psi$ has the transversality property.

\begin{example} \label{ex15}
{\rm Let $g$ be given by
$$
\begin{array}{cccc}
g \ : & \! \mathbb{R}^4  & \! \longrightarrow
& \! \mathbb{R}^{2}\\
& \! (x,y,z,w) & \! \longmapsto
& \! (x^2+y^2+z^3+zw^2,w^3+wz^2).
\end{array}
$$

The critical set of $g$ is $\{(x,y,0,0); \ x,y \in \mathbb{R}\}$ hence the discriminant is \linebreak $\{(x^2+y^2,0); \ x,y \in \mathbb{R}\}$ and the critical valued not isolated. Note that
$$V := g^{-1}(0) = \{(x,y,z,0) \in \mathbb{R}^4; \ x^2+y^2+z^3=0\}.$$

We find the points such that the fibers have tangency with $\mathbb{S}_\varepsilon^{3}$. Consider the matrix $M$ formed by the gradient vectors of the coordinate functions and the position vector:
$$ M = \left[
\begin{array}{cccccc}
2x & 2y & 3z^2+w^2 & 2wz \\
0  & 0  & 2zw & 3w^2+z^2 \\
x & y & z & w\\
\end{array}
\right] \,.$$
The computation of the minors
$$ \left|
\begin{array}{cccccccccc}
2u & 3z^2 + w^2  & 2wz\\
0& 2wz &3w^2+z^2\\
u& z & w\\
\end{array}
\right| = u(z^2+w^2)(3(z^2 +w^2)-2z),$$
with $ u = x, y$ we get,
\begin{equation}\label{eqx1}
\mbox{For}  \ u = x \ \mbox{we have}  \ \ x(z^2+w^2)(3(z^2 +w^2)-2z)
\end{equation}
\begin{equation}\label{eqy1}
\mbox{For}  \ u = y \ \mbox{we have}  \ \ y(z^2+w^2)(3(z^2 +w^2)-2z)
\end{equation}

We will develop the calculations for equation (\ref{eqy1}) and the procedure for equation (\ref{eqx1})  follows analogously.

The minor $y(z^2+w^2)(3(z^2 +w^2)-2z)$ is zero even if $y \neq 0$, so in general you have the points of the form
\begin{equation}\label{eqanalise}
(x,y,z,w) \ \ \mbox{such that} \ \ 3(z^2 +w^2)-2z=0.
\end{equation}

Let us see how they are. We can rewrite (\ref{eqanalise}) as $3z^2-2z+3w^2=0$ and see it as a quadratic equation in $z$, so we have
\begin{equation}\label{eqanalise1}
z=\dfrac{1 \pm \sqrt{1-9w^2}}{3}.
\end{equation}
Here $z$ is a real variable, se we need to have $1-9w^2 \geq 0$, which implies that $\frac{-1}{3} \leq w \leq \frac{1}{3}$. Set
$$w=\dfrac{t}{3}, \ \ \ \ \ \ t \in [-1,1]$$
and
$$z=\dfrac{1 \pm \sqrt{1-t^2}}{3}.$$
So these points have form
\begin{equation}\label{eqanalise2}
p(t)=\left(x,y,\dfrac{1 \pm \sqrt{1-t^2}}{3},\dfrac{t}{3}\right), \ \ \ \ \ \ t \in [-1,1],
\end{equation}
with $x$ and $y$ arbitrary.

So we have to analyze three cases:

\begin{enumerate}
  \item[\textbf{Case 1}] $z^2 +w^2 =0$ (this is, $z=0$ and $w=0$), then the points of the form $(x,y,0,0)$, which are critical points.
  \item[\textbf{Case 2}] $z^2 +w^2 \neq 0$ (this is, $z \neq 0$ or $w \neq 0$) and they satisfy equation (\ref{eqanalise}) (notice that in this case $z = 0$ implies $w = 0$ but not the other way
around). By (\ref{eqanalise1}) we know that, in this case we have points of the form (\ref{eqanalise2}). We need to prove that this kind of points cannot approach to $V$.
Note that taking in (\ref{eqanalise2}) the positive sign in the expression for $z$ we
have
$$p(1)=\left(x,y,\frac{1}{3}, \pm \frac{1}{3}\right) \ \ \mbox{and} \ \ p(0)=\left(x,y,\frac{2}{3},0\right),$$
so in principle we can have a curve of such points on a sphere taking the appropriate $x$ and $y$, for instance take
$$p(t)=\left(\dfrac{1-\sqrt{1-t^2}}{3},\dfrac{t}{3}, \dfrac{1 + \sqrt{1-t^2}}{3}, \dfrac{t}{3}\right), \ \ t \in [0,1],$$
and we have
$$||p(t)||^2= \dfrac{1-2 \sqrt{1-t^2} + 1 -t^2}{9}+\dfrac{t^2}{9}+\dfrac{1+2\sqrt{1-t^2}+1 - t^2}{9}+\dfrac{t^2}{9}=\dfrac{4}{9},$$

so the curve $p(t)$ with $t \in [0,1]$ lies on the sphere of radius $\frac{2}{3}$. If we
start at $p(1)$ and finish in $p(0)$, since $p(0)$ has its last coordinate $0$ in
principle it could be in $V$. But to be in $V$ it has to satisfy
\begin{equation}\label{eqanalise3}
x^2 + y^2 + z^3 = 0.
\end{equation}
Since $x^2 + y^2$ is always positive (unless $x = y = 0$ which is not the
case), in order that the point $(x,y,z,0)$ satisfies (\ref{eqanalise3}) it is necessary
to have $z^3$ negative, but in the case of the path $p(t)$ the $z$-coordinate
is always positive and therefore its cube $z^3$ is also always positive.
Hence this kind of points cannot approach to $V$.
  \item[\textbf{Case 3}] You have $z$ and $w$ such that $3(z^2 + w^2) - 2z \neq 0$ (do not satisfy
(\ref{eqanalise})). Then you must have $x = 0$ and $y = 0$ and you get points of the
form $(0,0,z,w)$, but when they approach to $V$  do not remain on a given sphere.
\end{enumerate}

Therefore it is not possible to have a sequence that satisfies the conditions in Theorem \ref{teopropdetransversalidade}, and then $g$ has the transversality property.

}
\end{example}

Another way to show that $g$ has the transversality property is to verify that $V \cap \Sigma_g = \{0\}$ by the following proposition.

\begin{proposition}\label{propdapropdetransversalidade}$($\rm{\cite{CMSS}}$)$ Let $f:(\mathbb{R}^n,0) \longrightarrow (\mathbb{R}^p,0), n > p \geq 1$, be a real analytic map-germ. If $V=f^{-1}(0) = \{0\}$ or if $\Sigma_f \cap V \subset \{0\}$, then $f$ has the transversality property.
\end{proposition}

\begin{example}\label{FMMPT}
{\rm Consider
$$
\begin{array}{cccc}
h \ : & \! \mathbb{R}^6  & \! \longrightarrow
& \! \mathbb{R}^{2}\\
& \! (x,y,z,w,t,r) & \! \longmapsto
& \! (x^2+y^2-z^2-w^2+t^3+tr^2,r^3+rt^2).
\end{array}
$$

The critical set of $h$ is: $\Sigma_{h}=\{(x,y,z,w,0,0); \ x,y,z,w \in \mathbb{R}\}$, hence \linebreak $ \Delta_{h} =\{(x^2+y^2-z^2-w^2,0); \ x,y,z,w \in \mathbb{R}\}$ so we have non-isolated critical values.

Note that $$V = h^{-1}(0) = \{(x,y,z,w,t,0); \ x^2+y^2-z^2-w^2+t^3=0\}\,.$$

Using  Proposition \ref{propdastangencias} we find the points such that the fibers have tangency with $\mathbb{S}_\varepsilon^{5}$. Consider the matrix:
$$ M = \left[
\begin{array}{cccccccccc}
2x & 2y & -2z & -2w & 3t^2+r^2 & 2tr \\
0  & 0  & 0 & 0 & 2tr & 3r^2+t^2 \\
x & y & z & w & t & r\\
\end{array}
\right].$$

The computation of the minors
$$ \left|
\begin{array}{cccccccccc}
2u & 3t^2 + r^2  & 2tr\\
0& 2tr &3r^2+t^2\\
u& t & r\\
\end{array}
\right|,$$
with $ u = x, y,-z,-w $ we get,
\begin{equation}\label{eqx}
\mbox{For}  \ u = x \ \mbox{we have}  \ \ x(t^2+r^2)(3(t^2 +r^2)-2t)
\end{equation}
\begin{equation}\label{eqy}
\mbox{For}  \ u = y \ \mbox{we have}  \ \ y(t^2+r^2)(3(t^2 +r^2)-2t)
\end{equation}
\begin{equation}\label{eqz}
\mbox{For}  \ u = -z \ \mbox{we have}  \ \ z(t^2+r^2)(3(t^2 +r^2)+2t)
\end{equation}
\begin{equation}\label{eqw}
\mbox{For}  \ u = -w \ \mbox{we have}  \ \ w(t^2+r^2)(3(t^2 +r^2)+2t)
\end{equation}

The obvious cases when these minors are zero is when either $x=y=z=w=0$ getting points of the form
\begin{equation}\label{eqs}
(0, 0, 0, 0, t, r)
\end{equation}
or $t = r =0$, getting points of the form
\begin{equation}\label{eqss}
(x,y,z,w,0,0).
\end{equation}

Now if $x^2 + y^2 \neq 0$ in order to have a point with $t^2 + r^2 \neq 0$, by (\ref{eqx}) and (\ref{eqy}), we need that $t$ and $r$ satisfy the equation $3(t^2 + r^2) - 2t$ as in Example \ref{ex15}. So we have points of the form
\begin{equation}\label{eqsx}
p(s) = \left( x,y,0,0, \frac{1 \pm \sqrt{1-s^2}}{3}, \frac{1s}{3} \right), \ s \in [-1,1].
\end{equation}

Analogously if $z^2 +w^2 \neq 0$ in order to have a point with $t^2 + r^2 \neq 0$ by (\ref{eqz}) and (\ref{eqw}), we need that $t$ and $r$ satisfy the equation $3(t^2 + r^2) + 2t$ (notice the difference on sign). So we have points of the form
\begin{equation}\label{eqsz}
q(s) = \left( 0,0,z,w, \frac{-1 \pm \sqrt{1-s^2}}{3}, \frac{1s}{3} \right), \ s \in [-1,1].
\end{equation}

Note that we cannot have points with $x^2 + y^2 \neq 0$ and $z^2 + w^2 \neq 0$ at the same time, because by (\ref{eqx}) and (\ref{eqy}) the coordinate $t$ is always positive, while by (\ref{eqz}) and (\ref{eqw}) the $t$ coordinate is always negative, so the only possibility for this is $t=0$ which implies $r=0$ but in this case we get $(x,y,z,w,0,0).$
It is easy to check that all the other minors vanish on the points of the forms (\ref{eqs}), (\ref{eqss}), (\ref{eqsx}) and (\ref{eqsz}). It is trivial to see that the points of the form (\ref{eqs}) cannot approach to $V$, the points of type (\ref{eqss}) are critical points. Only remains to check the other two types.

As in Example \ref{ex15}, when $s \to 0$ we have $p(s) \to (x,y,0,0, \frac{2}{3},0)$ and this point is in $V$ if it satisfies $x^2 + y^2 + t^3 =0$, but since $x^2+y^2 >0$ we need $t < 0$ but in $p(s)$ the t-coordinate is always positive. Analogously, when $s \to 0$ we have $q(s) \to (0,0,z,w,\frac{-2}{3},0)$ and this point is in $V$ if it satisfies $-z^2 - w^2 + t^3=0$, but since $-z^2 - w^2 < 0$ we need $t > 0$ but in $q(s)$ the t-coordinate is always negative.

Note that $ h $ has the transversality property and also in this example \linebreak $V \cap \Sigma_{h} \neq \{0\}$, because the point $(1,1,1,1,0,0) \in V \cap \Sigma_{h}$. Thus in this case it is not possible to prove the transversality property using Proposition \ref{propdapropdetransversalidade}.
}
\end{example}

\section{The main theorem}
\label{sec:3}

Consider a map $\psi: \mathbb{C}^{n} \longrightarrow \mathbb{C}$ of the form
$$\psi(z) = \sum_{j=1}^n \lambda_j  z^{a_j}_{j}\overline{z_j}^{b_j}, \ \ a_j,b_j \in \mathbb{Z}, a_j,b_j \geq 0,$$
with $a_j = b_j$ for at least one index $j$ and $a_j \neq b_j$ for
at least one index $j$ for $j = 1,\cdots,n$ and $\lambda_j \in \mathbb{C}^*$, $\lambda_j = \alpha_j + i \beta_j, \ \alpha_j,\beta_j \in \mathbb{R}$.
We want to determine the conditions for $\psi$ to have Milnor Fibration.

Note that it is easy to see that in the case $a_j=0,b_j=1$ or $a_j=1,b_j=0$ with $Re(\lambda_j)\neq 0$ and $Im(\lambda_j) \neq 0$, $j=1, \cdots, n$, $\psi$ is a submersion, and therefore we do not need to investigate these situation.

We know from \cite{Ruas-Seade-Verjovsky} that if $a_j \ne b_j$ for all $j = 1, \cdots, n$, then these maps have an isolated critical point.

The above polynomial $\psi(z)$ are not necessarily holomorphic maps, so \linebreak we consider then as real analytic maps $\psi: \mathbb{R}^{2n} \to \mathbb{R}^2$, where \linebreak $z_j = x_j + i y_j, \ x_j,y_j \in \mathbb{R}$, $j=1, \cdots, n$.  We set $\psi_1 (z) = Re(\psi(z))$ and $\psi_2 (z) = Im(\psi(z))$.

We say that an index $j \in \{1,2, \cdots, n\}$ is a critical index of $\psi$ if $a_j = b_j$. (Because these
are the indices which give critical points). Let $C \subset \{1,...,n\}$ be the subset of critical indices of $\psi$.
Let $j,k \in \{1,...,n\}$, we say that $\lambda_j$ and $\lambda_k$ are colinear if $\alpha_j \beta_k = \beta_j \alpha_k$ (*).

In fact we can see this geometrically. We have $\lambda_j = \alpha_j + i \beta_j$ for $j = 1,2$ and
$\dfrac{\beta_j}{\alpha_j}= \arctan \lambda_j$ is the tangent of the argument of $\lambda_j$. If equation (*) is satisfied, this implies
$$\arctan \lambda_1 = \dfrac{\beta_1}{\alpha_1} =\dfrac{\beta_2}{\alpha_2}= \arctan \lambda_2$$
which means that either
$$ arg \lambda_2 = arg \lambda_1$$
or
$$arg \lambda_2 = arg \lambda_1 + \pi.$$

In any case, $\lambda_1$ and $\lambda_2$ are on the same (real) line, this is, are collinear .


\begin{lemma} Colinearity of critical indices is an equivalence relation on the
set $C$ of critical indices of $\psi$.
\end{lemma}
\vspace{-0.2cm}
\begin{prova} The collinearity of the critical indices is clearly reflective and \linebreak symmetrical so we will prove only the transitivity relation.

Consider that $a_j = b_j$ for $j = 1,2,3$ and set
$\lambda_j = \alpha_j + i \beta_j$ for $j = 1,2,3$ and suppose that
\begin{equation}\label{ex12}
\alpha_1 \beta_2 = \beta_1 \alpha_2 \ \ then \ \ \lambda_1 \ and \ \lambda_2 \ are \ collinear
\end{equation}
\begin{equation}\label{ex13}
\alpha_1 \beta_3 = \beta_1 \alpha_3 \ \ then \ \ \lambda_1 \ and \ \lambda_3 \ are \ collinear
\end{equation}
and we want to show that
\begin{equation}\label{ex23}
\alpha_2 \beta_3 = \beta_2 \alpha_3 \ \ then \ \ \lambda_2 \ and \ \lambda_3 \ are \ collinear
\end{equation}

Notice that the first two conditions (equations (\ref{ex12}) and (\ref{ex13})
above) imply equation (\ref{ex23}): multiply (\ref{ex12}) and (\ref{ex13}) to get
$$\alpha_1 \beta_2 \beta_1 \alpha_3 = \beta_1 \alpha_2 \alpha_1 \beta_3,$$
and canceling $\alpha_1 \beta_1$ in both sides we get (\ref{ex23}). So if two pairs of indices
$\{k,l\}$ and $\{k,m\}$ satisfy $\alpha_k \beta_l = \beta_k \alpha_l$ and $\alpha_k \beta_m = \beta_k \alpha_m$ and both pairs have an index in
common, then the pair $\{l,m\}$ also satisfies the condition.

Hence, one can define an equivalence relation on the set of indices $j$ such that $a_j = b_j$ saying
that two indices are related if they are in the same line.
\end{prova}

\vspace{0.3cm}

Denote by $\mathcal{C}$ the set of equivalence classes of colinear critical indices. Let
$J \in \mathcal{C}$ be an equivalence class, we say that $J$ has size $r$ if the cardinality
$|J|$ of $J$ is $r$.

Given an equivalence class of critical indices $J$ of size $r$ define the set
$$\Sigma_J = \{(0,\cdots,z_{j_1}, \cdots,z_{j_2},\cdots,z_{j_r}, \cdots,0) \in \mathbb{C}^n; j_l \in J \ \mbox{for} \ l = 1,\cdots,r\}.$$

Note that $\Sigma_J$ is the subspace of $\mathbb{C}^n$ of dimension $r$ generated by the coordinates $z_j$ with $j \in J$.

\begin{proposition}\label{proppontoscriticos}
Consider
$$\psi(z)= \sum_{j=1}^n \lambda_{j} z^{a_j}_{j}\overline{z_j}^{b_j}, \ a_j,b_j \in \mathbb{Z}, a_j,b_j \geq 0,$$
and with $a_j = b_j$ for at least one index $j$ and $a_j \neq b_j$ for
at least one index $j$ for $j = 1,\cdots,n$, that is, $0 < |C| < n$ where $C$ is the
set of critical indices. Then the critical set $\Sigma_\psi$ of $\psi$ is given by
$$\Sigma_\psi = \bigcup_{J \in C} \Sigma_J.$$
%
%

\end{proposition}

\vspace{-0.2cm}
\begin{prova}
Let  $\psi_1 (z) = Re(\psi(z))$ and $\psi_2(z) = Im(\psi(z))$. Assuming we have only a single index with equal exponents (only critical index of $\psi$, $J=\{j\}$), in this case $z_j^{a_j}\overline{z_j}^{a_j} = ||z_j||^{2 a_j}$, and calculating the jacobian matrix of this situation we conclude directly that $\Sigma_{\psi} = \{(0, \cdots, z_j, \cdots, 0)\}$  .

Without loss of generality, assume $a_j=b_j$ for $j=1,2, \cdots, r$ (this is, $j$ is a critical index for $j = 1,2 ,\cdots,r)$. Thus we have
$$\psi(z) = \sum_{j=1}^r \lambda_j z^{a_j}\overline{z_j}^{a_j}+ \sum_{j=r+1}^n \lambda_{j} z^{a_j}\overline{z_j}^{b_{j}}.$$

Setting $\lambda_{j} = \alpha_{j} + i \beta_{j}, \ \alpha_j, \beta_j \in \mathbb{R}^*$, for $j=1,\cdots, r$, we get
$$\psi_1(z) = \sum_{j=1}^r \alpha_{j} ||z_{j}||^{2a_{j}} + \dfrac{1}{2} \sum_{j=r+1}^n (\lambda_j z_j^{a_j}\overline{z_j}^{b_j} + \overline{\lambda_j}\overline{z_j}^{a_j}z_j^{b_j})$$
and
$$\psi_2(z) =  \sum_{j=1}^r \beta_{j} ||z_{j}||^{2a_{j}}  + \dfrac{1}{2} \sum_{j=r+1}^n (\lambda_j z_j^{a_j}\overline{z_j}^{b_j} - \overline{\lambda_j} \overline{z_j}^{a_j}z_j^{b_j}).$$

As discussed in \cite{Ruas-Seade-Verjovsky}, let us consider $ 2 \times 2$ minor of $\psi$ corresponding to the partial derivatives of $(\psi_1,\psi_2)$ with respect to $z_i$ and $\overline{z_i}$. Its determinant is:
$$||\lambda_i ||^2a_{i}^2||z_{i}||^{2a_{i}-2}z_{i} ||z_{i}||^{2b_i}- ||\lambda_i ||^2b_{i}^2||z_{i}||^{2a_{i}}z_{i} ||z_{i}||^{2b_i-2},$$
so this is $0$ iff $z_i=0$ or $a_i=b_i$.

Suppose $j_1, j_2$, without loss of generality, are the indeces such that we have that the exponents are equal (this is, $j_1$ and $j_2$ are a critical index). Now let's check the jacobian matrices 2 $\times$ 2 of $\psi$ in relation to those partial derivatives from $z_j$ and $\overline{z_j}$, we have
$$a_{j_1}a_{j_2}||z_{j_1}||^{2a_{j_1}-2}z_{j_1} ||z_{j_2}||^{2a_{j_2}-2}\overline{z_{j_2}}(\alpha_{j_1}\beta_{j_2} - \alpha_{j_2}\beta_{j_1})$$
and
$$a_{j_1}a_{j_2}||z_{j_1}||^{2a_{j_1}-2}\overline{z_{j_1}} ||z_{j_2}||^{2a_{j_2}-2}z_{j_2}(\alpha_{j_1}\beta_{j_2} - \alpha_{j_2}\beta_{j_1}),$$
are zero when $z_j$'s are zero or $(\alpha_{j_1}\beta_{j_2} - \alpha_{j_2}\beta_{j_1}) = 0$, when second statement \linebreak happens then $\Sigma_\psi = \{(0, \cdots,z_{j_1}, 0,\cdots,z_{j_2}, \cdots, 0)\}$. When $(\alpha_{j_1}\beta_{j_2} - \alpha_{j_2}\beta_{j_1}) \neq 0$ we have $\Sigma_\psi = \{(0, \cdots, z_{j_1},\cdots,0)\} \cup \{(0,\cdots, z_{j_2},\cdots, 0)\}$, proving the lemma.



\end{prova}

\vspace{0.3cm}

\begin{lemma}\label{lemapsiradial}
The map $\psi$ is radial weighted homogeneous
of type $(p1,...,pn;a)$ with $a =  lcm(a_1+b_1, \cdots, a_n+b_n)$ and $p_j = \frac{a}{a_j + b_j} > 0$.
\end{lemma}

\begin{prova} Using the $\mathbb{R}^+$-action defined in \cite{CA}, given by
$$t \cdot (z_1, \cdots, z_n) = (t^{p_1}z_1, \cdots, t^{p_n}z_n), \ t \in \mathbb{R}^+,$$
where $P = (p_1, \cdots, p_n)$ are the radial weights. It is weighted homogeneous of \linebreak
degree $a$ such that $\psi(t \cdot (z)) = t^{a} \psi(z), \ t \in \mathbb{R}^+$. We define \linebreak $a =  lcm(a_1+b_1, \cdots, a_n+b_n)$ and $p_j = \frac{a}{a_j + b_j} > 0$.
\end{prova}

\begin{definition}\label{discriminantelinear}$($\rm{\cite{CMSS}}$)$ {\rm A map-germ $f:(\mathbb{R}^n,0) \longrightarrow (\mathbb{R}^p,0)$ of class $\mathcal{C}^\ell$ has} {linear discriminant} {\rm if there exists $\eta= \eta(f) > 0$ such that the intersection of the discriminant $\Delta_f$ of $f$ with the ball $\mathbb{B}_\eta^p$ is a union of line-segments, that is:
$$\Delta_f \cap \mathbb{B}_\eta^p = \mbox{Cone}(\Delta_f \cap \mathbb{S}_\eta^{p-1}).$$
We call $\eta$ a \emph{linearity radius} of $\Delta_f$.}
\end{definition}

The next result shows that $\psi$ has linear discriminant.

\begin{proposition}\label{caracdodiscriminante}
Let $$\psi(z)=  \sum_{j=1}^n \lambda_j z^{a_j}_j\overline{z_j}^{b_j}, \ \ a_j,b_j \in \mathbb{Z}, a_j,b_j \geq 0, \lambda_j \in \mathbb{C}^*.$$
Let $C$ be the set of critical indices of $\psi$ and suppose $0 < |C| < n$. For
each $j \in C$ let $\mathcal{L}_j$ be the real ray in $C$ from $0$ to $\lambda_j$. Then we have:
\begin{enumerate}
  \item[1] The discriminant $\Delta_\psi$ of $\psi$ is given by
  $$\Delta_\psi = \bigcup_{j \in C} \mathcal{L}_j.$$
  If several $\lambda_j$ (say $r$ of them) are on the same line, their corresponding indices form an equivalence class $J$ of critical colinear indices of size $r$.
  \item[2] There are no complete lines in $\Delta_\psi$, if and only if for every equivalence class $J$ of critical colinear indices, all the indices $j \in J$ have the same argument.
\end{enumerate}
\end{proposition}
\vspace{-0.2cm}
\begin{prova} Let $J = \{j_1, \cdots,j_r\}$ be an equivalence class of critical colinear
indices of size $r$. Then we have that $\lambda_{j_1}, \cdots, \lambda_{j_r}$ are colinear in $C$. Set $\theta = \arg \lambda_{j_1}$, so, $\lambda_{j_l} = \mu_{j_l} e^{i\theta}$ with $\mu_{j_l} \in \mathbb{R}$ for $l = 1,\cdots,r$. Evaluating $\psi$ in a critical point
$$p = (0,\cdots,z_{j_1},\cdots,z_{j_2},\cdots,z_{j_r}, \cdots,0) \in \Sigma_J,$$
we obtain
\begin{equation}\label{discriminanteeq}
\psi(p)= \sum_{l=1}^r \lambda_{j_l} ||z_{j_l}||^{2a_{j_l}} = \left(\sum_{l=1}^r \mu_{j_l}||z_{j_l}||^{2a_{j_l}}\right)e^{i\theta}.
\end{equation}

Since $(\sum_{l=1}^r \mu_{j_l}||z_{j_l}||^{2a_{j_l}})$ is a real number, as $p$ varies in $\Sigma_J$ we obtain
the line (or ray if all $\mu_{j_l} > 0$) which passes through $0$ and $e^{i \theta}$. Each
equivalence class of critical colinear indices contributes with one line (or
ray) of the discriminat $\Delta_\psi$: if two equivalence classes $J_1$ and $J_2$ contribute
the same line, the definition of the colinearity equivalence relation implies
that $J_1 = J_2$.

For the last claim of the Proposition \ref{proppontoscriticos}, if $J = \{j_1,\cdots,j_r\}$ is an equivalence
class of critical colinear indices, all $\lambda_{j_1},\cdots, \lambda_{j_r}$ have the same argument, if
and only if $\mu_{j_l} > 0$ for $l = 1,\cdots,r$ and $(\sum_{l=1}^r \mu_{j_l}||z_{j_l}||^{2a_{j_l}})$ in equation (\ref{discriminanteeq})
is a positive real number, so as $p \in \Sigma_J$ varies, we only obtain the ray from $0$ through $\lambda_1$.
\end{prova}

\vspace{0.2cm}


Given an analytic map-germ as above,  $f: (\mathbb R^n,\0) \to (\mathbb R^p,0)$, $n > p \ge 1$, with a critical point at the origin $\0 \in \mathbb R^n$ and such that $V:= f^{-1}(0)$ has dimension greater then zero. We now consider fibrations on small spheres. One can associate to $f$ a canonical pencil as in \cite{CSS}: to each line $\mathcal{L}$ through the origin in $\mathbb R^p$ we associate the analytic set
$$X_{\mathcal{L}} := \{x \in \mathbb R^n \, \big | \, f(x) \in \mathcal{L} \, \}\;.
$$

The union of all these analytic sets is the whole ambient space, their intersection is $V$ and each of these is non-singular away from the critical set of $f$. We now consider map-germs with non-isolated critical values but with linear discriminant. The concept of regularity was extended to this setting as follows:

\begin{definition}\label{definicaod-regular} {\rm (\cite{CMSS})} Let $f:(\mathbb{R}^n,0) \longrightarrow (\mathbb{R}^p,0)$ be a map-germ of class $\mathcal{C}^\ell$ with $\ell \geq 1$ with linear discriminant. We say that $f$ is {$d-$regular} at $0$ if for any representative $f$  there exists $\varepsilon_{0}>0$ small enough such that $f(\mathbb{B}_{{\varepsilon}_{0}}^{n})\subset \stackrel{\circ}{\mathbb{B}^{p}_{\eta}}$,  where $\eta$ is a linearity radius for $f$, such that  for each $X_{\mathcal{L}} \setminus f^{-1}(\Delta_f)$. There exists $\varepsilon_0 >0$ such that every sphere of radius $\varepsilon \le \varepsilon_0$ and center at $\0$ intersects  transversally each $X_{\mathcal{L}} \setminus f^{-1}(\Delta_f)$.
\end{definition}

\begin{theorem}\label{teodafibracaoCMSS}$($\rm{\cite{CMSS}}$)$ Let $f:(\mathbb{R}^n,0) \longrightarrow (\mathbb{R}^p,0)$ with $n \geq p \geq 2$ be a map-germ of class $\mathcal{C}^\ell$, $\ell \ge 1$, with linear discriminant and linearity radius $\eta > 0$. If $f$ is $d$-regular and has the transversality property, then there exists a positive real number $0 < \eta \ll \varepsilon$ such that:
$$\phi := \frac{\psi}{||\psi||} : \mathbb{S}_\varepsilon^{n-1} \setminus W \longrightarrow \mathbb{S}^{p-1} \setminus \mathcal{A},$$
is a (differentiable) locally trivial fibration over its image, where $W:= f^{-1}(\Delta_f^\eta)$, $\mathcal{A}_\eta = \Delta_\psi \cap \mathbb{S}_\eta^{2n-1}, \pi(\mathcal{A}_\eta) =\mathcal{A}$, with $\pi$ the projection onto the unit sphere. If $f$ is real analytic then the fibration $\phi$ is smooth indeed. Moreover, this fibration is equivalent to the local Milnor-L\^e Fibration.

\end{theorem}

\begin{proposition}
\label{propfibracaodafamiliad-regular}

The map $$\psi(z)= \sum_{j=1}^n \lambda_j  z^{a_j}_{j}\overline{z_j}^{b_j}, \ \ a_j,b_j \in \mathbb{Z}, a_j,b_j \geq 0,$$
with $a_j = b_j$ for at least one $j$ and $a_j \neq b_j$ for at least one $j$, $j = 1, \cdots, n$, is $d$-regular.
\end{proposition}
\vspace{-0.2cm}
\begin{prova} It follows from the fact that $\psi$ has $\mathbb{R}^+$-action (Lemma \ref{lemapsiradial}) and doing the same procedure developed in the proof of Proposition 3.4 in \cite{Cisneros}.
\end{prova}

\vspace{0.5cm}

By Proposition \ref{propdapropdetransversalidade} in the next result we present conditions for $\psi$ to have the transversality property .

\begin{lemma} \label{lemad-regular}
Let
$$\psi(z)= \sum_{j=1}^n \lambda_j z^{a_j}_{j}\overline{z_j}^{b_j}, \, a_j,b_j \in \mathbb{Z}, a_j,b_j > 0, \lambda_j \in \mathbb{C}^*.$$

Let $C$ be the set of critical indices of $\psi$ and suppose $0 < |C| < n$. Then $\psi$ has the transversality property if it satisfies the following (equivalent) conditions:
\begin{enumerate}
 \item[i)] The discriminant $\Delta_\psi$ does not have any complete line.
 \item[ii)] For every equivalence class $J$ of critical colinear indices, all the indices $j \in J$ have the same argument.
\end{enumerate}

\end{lemma}
\vspace{-0.2cm}
\begin{prova} The equivalence of the two conditions is established in Proposition \ref{caracdodiscriminante}, using condition $2$ we prove that $\psi$ satisfies the criterion of Proposition \ref{propdapropdetransversalidade}, \linebreak this is, that $\Sigma_\psi \cap V = \{0\}$, so we need to check that if $p$ is a critical point of $\psi$ then $\psi(p) \neq 0$ unless $p = 0$. Let $p$ be a critical point of $\psi$, by Proposition \ref{proppontoscriticos}, $p \in \Sigma_J$ for some equivalence class of critical colinear indices $J = \{j_1, \cdots, j_r\}$. By condition $2$ (Proposition \ref{caracdodiscriminante}), all the $\lambda_{j_l}$ for $j_l \in J$ have the same argument, thus, if $\theta= \arg \lambda_{j_l}$, we have that $\lambda_{j_l} = \mu_{j_l}e^{i \theta}$ with $0 < \mu_{j_l} \in \mathbb{R}$ for $l=1, \cdots, r$. By equation (\ref{discriminanteeq}) we have that $\psi(p) = (\sum_{l=1}^{r} \mu_{j_l}||z_{j_l}||^{2_{a_{j_l}}})e^{i \theta}$ and $(\sum_{l=1}^{r} \mu_{j_l}||z_{j_l}||^{2_{a_{j_l}}})$ is a positive real number. Then $\psi(p) = 0$ implies that $z_{j_l} = 0$ for every $j_l \in J$, but since $p \in \Sigma_J$ we also have that $z_j = 0$ for every $j \in\!\!\!\!\!/ \ J$, so $p = 0$.

Hence, $\Sigma_\psi \cap V = \{0\}$ and by Proposition \ref{propdapropdetransversalidade} $\psi$ has the transversality \linebreak property.
\end{prova}

\begin{theorem}\label{teofibracaoitens}
The map
$$\phi = \frac{\psi}{||\psi||}:(\mathbb{S}_\varepsilon^{2n-1} \setminus \psi^{-1}(\Delta_\psi^\eta) ) \longrightarrow\mathbb{S}^1 \setminus \mathcal{A},$$
where $\mathcal{A}_\eta = \Delta_\psi \cap \mathbb{S}_\eta^{2n-1}, \pi(\mathcal{A}_\eta) =\mathcal{A}$, with $\pi$ the projection onto the unit sphere, and $\Delta_\psi^n = \Delta_\psi \cap \mathbb{B}^{2n-1}_\eta$,is a locally trivial fibration over its image under the following conditions:
\begin{enumerate}
 \item[i)] The discriminant $\Delta_\psi$ does not have any complete line.
 \item[ii)] For every equivalence class $J$ of critical colinear indices, all the indices $j \in J$ have the same argument.
\end{enumerate}\end{theorem}
\vspace{-0.2cm}
\begin{prova} The proof follows from Proposition \ref{propfibracaodafamiliad-regular} and Lemma \ref{lemad-regular} by Theorem 3.8 in \cite{CMSS}.

\end{prova}

The next example illustrates a map that satisfies the last theorem.

\begin{example}
Let $f(z)=(1+i) z_1 \overline{z_1}+ (-2-i) z_2^2\overline{z_2}^2+ i z_3^2 \overline{z_3}$. Taking \linebreak $z_j = x_j + i y_j, \ x_j , y_j \in \mathbb{R}, \ j=1,2,3$ one has
$$f_1(x_1,y_1,x_2,y_2,x_3,y_3)=x_1^2+y_1^2-2(x_2^2+y_2^2)^2-x_3^2y_3-y_3^3$$
and $$f_2(x_1,y_1,x_2,y_2,x_3,y_3)=x_1^2+y_1^2 - (x_2^2+y_2^2)^2 + x_3^3+x_3y_3^2$$
and the jacobian matrix of $f$ corresponding to the partial derivatives of $f_1$ and $f_2$
$$ Jf = \left[
\begin{array}{cccccc}
2x_1 & 2y_1 & -8(x_2^2+y_2^2)x_2 & -8(x_2^2+y_2^2)y_2 & -2x_3y_3 & -x_3^2-3y_3^2\\
2x_1 & 2y_1 & -4(x_2^2+y_2^2)x_2 & -4(x_2^2+y_2^2)y_2 & 3x_3^2+y_3^2 & 2x_3y_3\\
\end{array}
\right]. $$

Hence $\Sigma_f=\{(x_1,y_1,0,0,0,0); x_1,y_1 \in \mathbb{R}\} \cup \{(0,0,x_2,y_2,0,0); x_2,y_2 \in \mathbb{R}\}$ and then
$$\Delta_f =\{(x_1^2+y_1^2, x_1^2+y_1^2); x_1,y_1 \in \mathbb{R}\} \cup \{(-2(x_2^2+y_2^2)^2, - (x_2^2+y_2^2)^2); \ x_2,y_2 \in \mathbb{R}\}$$
that is a union of to line-segments. Then by condition i) of Lemma \ref{lemad-regular} $f$ has the transversality property. Therefore $f$ satisfies the conditions of Theorem \ref{teofibracaoitens} and thus we have that $f$ has the Milnor Fibration.

\end{example}

In \cite{Cisneros-Aurelio} and \cite{CMSS} the authors prove a result that ensures equivalence between the Milnor-L\^{e} Fibration on the tube and the Milnor Fibration on the sphere in this context. The next proposition is a consequence of this result.

\begin{proposition}\label{teoequivalencia}
Suppose that $\psi$ satisfies the Theorem \ref{teofibracaoitens}. The fibrations
$$\tilde{\psi}:= \tilde{\pi} \circ \psi_{|}: \mathbb{B}_\varepsilon^{2n} \cap \psi^{-1}(\mathbb{S}^1_\delta\setminus \Delta_\psi) \longrightarrow \mathbb{S}^1 \setminus \mathcal{A},$$
and
$$\phi :\mathbb{S}_\varepsilon^{2n-1} \setminus W \longrightarrow\mathbb{S}^1 \setminus \mathcal{A},$$
are equivalent, where $\tilde{\pi}: \mathbb{B}^{p}_{\delta} \to \mathbb{S}^{p-1}$ is the radial projection.
\end{proposition}

We show that for $\psi$, mixed homogeneous polynomials the Milnor-Lê fiber on a critical value is homeomorphic to a fiber on a regular value, \cite{AR}.

We know of the existence of fibration in the tube for critical values, due to the following result:

Consider $X \subset \mathbb{R}^m$ and $Y \subset \mathbb{R}^n$ subanalytic sets, with $m \geq n$. Let
$$f :X \to Y$$
be a subanalytic map that extends to a continuous map $\tilde{f}: \mathbb{R}^m \to \mathbb{R}^n$.

We have:

\begin{proposition}$($\rm{\cite{AR}}$)$
Let $B \subset \mathbb{R}^m$ be a compact subanalytic set and let $\mathcal{C} \subset \mathbb{R}^n$ be a one-dimensional subanalytic set contained in $f(X \cap B)$. For every $y \in \mathcal{C}$ there exists a positive real number $\delta>0$ such that the restriction:
$$f_|: f^{-1} \big( (\mathcal{C} \backslash \{y\}) \cap \mathbb{D}_\delta(y) \big) \cap B \to (\mathcal{C} \backslash \{y\}) \cap \mathbb{D}_\delta(y)$$
is the projection of a topological trivial fibration, where $\mathbb{D}_\delta(y)$ denotes the closed ball of radius $\delta$ around $y$ in $\mathbb{R}^n$.
\end{proposition}

\vspace{0.2cm}

In the next example, one can see that the manifolds defined by $\psi^{-1}(d), \ d \in \mathbb{R}^2$ and $d\neq (0,0)$, are all smooth manifolds.

The fiber $f^{-1}(0,0)$ is singular, but it is possible to show that $f^{-1}(0,0)$ is homeomorphic to $f^{-1}(d), d\neq (0,0)$. It means 
that the topological information is preserved even in the special fiber. 

\begin{example} Consider $f(z)= z_1\overline{z_1}+z_2^2\overline{z_2}$ in $\mathbb{R}^2$, we take $f_1(z)= Re(f(z))$, $f_2(z) = Im(f(z))$ and $z_j = x_j + i y_j, \ x_j , y_j \in \mathbb{R}, \ j=1,2$,
$$f_1(x_1,y_1,x_2,y_2) = x_1^2+y_1^2+x_2(x_2^2+y_2^2) \ \ \ \mbox{and} \ \ \ f_2(x_1,y_1,x_2,y_2)=y_2(x_2^2+y_2^2)$$

In this situation $\Sigma_f= \{(x_1,y_1,0,0); \ x_1, y_1 \in \mathbb{R}\}$ and consequently $\Delta_f= \linebreak \{(x_1^2+y_1^2,0); \ x_1,y_1 \in \mathbb{R}\} = \{(c,0), c \in \mathbb{R}^+\}$.

Recall that the discriminant of $f$ has fibration.

Let us take $c \in \mathbb{R}^+$, in the situation we have
\begin{align*}
f^{-1}(c,0) & = \{(x_1,y_1,x_2,y_2) \in \mathbb{R}^4; \ x_1^2+y_1^2+x_2(x_2^2+y_2^2) = c \ \mbox{and} \ y_2(x_2^2+y_2^2)=0\} \\
& = \{(x_1,y_1,x_2,0) \in \mathbb{R}^4; \ x_1^2+y_1^2+x_2^3=c\},\\
f^{-1}(0,0) & = \{(x_1,y_1,x_2,0) \in \mathbb{R}^4; \ x_1^2+y_1^2+x_2^3=0\} \ \ \  and \\
f^{-1}(-c,0) & = \{(x_1,y_1,x_2,0) \in \mathbb{R}^4; \ x_1^2+y_1^2+x_2^3=-c\}
\end{align*}

The first fiber of the discriminants, the second the special fiber and finally the fiber of a regular value opposite the critical value used. Through the image below, considering $ \mathbb{R}^4$ a cut $y_2= 0$, you can see a representation of the three fibers described above.

\begin{figure}[h]
\begin{center}
\includegraphics[scale=0.35,angle=0]{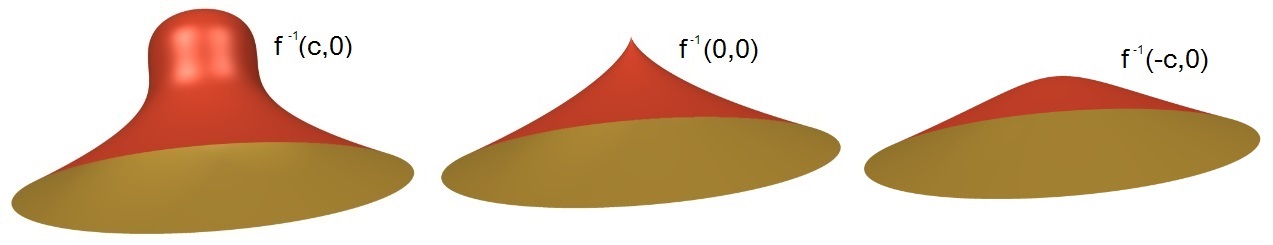}
 \caption{The three homeomorphic fibers}
  \end{center}
  \end{figure}

\end{example}

\section{Special Case}

Note that the Example \ref{FMMPT} is a this is case not covered by Theorem \ref{teofibracaoitens}. The map in that example is given by,
$$h(x,y,z,w,t,r) = \psi(z)=z_1 \overline{z_1} - z_2 \overline{z_2} + z_3^2\overline{z_3},$$
where $z_1 = x+iy, \ z_2 = z+i w, \ z_3 = t+ir, \ x,y,z,w,t,r \in \mathbb{R}$.

For a class of $\psi$ containg example above, we proof the existence of Milnor Fibration using the d-regularity (Proposition \ref{propfibracaodafamiliad-regular}) and proving by Proposition \ref{teopropdetransversalidade} that we have in this case the transversality property.

\begin{lemma}\label{lemacasoespecial}
Consider mixed polynomials given by
\begin{equation}\label{eq21}
\psi(z)=\sum_{j=1}^{r}\mu_j z_j^{a_j}\overline{z_j}^{a_j} - \sum_{j=r+1}^{n-1}\mu_j z_j^{a_j}\overline{z_j}^{a_j} \pm \mu_n z_n^2 \overline{z_n}
\end{equation}
or
\begin{equation}\label{eq12}
\psi(z)=\sum_{j=1}^{r}\mu_j z_j^{a_j}\overline{z_j}^{a_j} - \sum_{j=r+1}^{n-1}\mu_j z_j^{a_j}\overline{z_j}^{a_j} \pm \mu_n z_n \overline{z_n}^2,
\end{equation}
with $a_j, b_j \in \mathbb{Z}, a_j, b_j > 0$ and $0 < \mu_j \in \mathbb{R}$ for $j = 1,\cdots,n$. Then the map
$$\phi = \frac{\psi}{||\psi||}:(\mathbb{S}_\varepsilon^{2n-1} \setminus W) \longrightarrow\mathbb{S}^1 \setminus \mathcal{A},$$
is a locally trivial fibration, where $\Delta_\psi^\eta= \Delta_\psi \cap \mathbb{B}^2_\eta, W= \psi^{-1}(\Delta_\psi^\eta), \mathcal{A}_\eta=\Delta_\psi^\eta \cap \mathbb{S}^1_\eta$ and $\mathcal{A} = \pi(\mathcal{A}_\eta)$ and $\pi:\mathbb{C}^* \to \mathbb{S}^1$ is the radial projection.
\end{lemma}
\vspace{-0.2cm}
\begin{prova} We give the proof for (\ref{eq21}) with positive $\mu_n z_n^2 \overline{z}_n$ term, the proof of the other case and of (\ref{eq12}) are analogous. Note that $j = 1,\cdots, n-1$ are the critical indices of $\psi$ and since $\mu_j \in \mathbb{R}$ for $j = 1,\cdots,n - 1$ all of them form the only equivalence class of colinear critical indices. Taking $z_j = x_j + iy_j, \ x_j, y_j \in \mathbb{R}$ the real and imaginary parts of $\psi$ are given by
$$\psi_1(x_1, \cdots,y_n)=\sum_{j=1}^{r}\mu_j (x_j^2+y_j^2)^{a_j} - \sum_{j=r+1}^{n-1}\mu_j (x_j^2+y_j^2)^{a_j} + \mu_n x_n(x_n^2 +y_n^2)$$
$$\psi_2(x_1,\cdots,y_n)=\mu_n y_n(x_n^2 + y_n^2).$$
We have that
\begin{align*}
V & = \psi^{-1}(0) \\
& = \left\{(x_1,\cdots,y_{n-1},x_n,0)|\sum_{j=1}^{r}\mu_j (x_j^2+y_j^2)^{a_j} - \sum_{j=r+1}^{n-1}\mu_j (x_j^2+y_j^2)^{a_j}+ \mu_n x_n^3=0\right\}.
\end{align*}
We want to use Proposition \ref{teopropdetransversalidade}, so we
need to compute for which points its fibre is tangent to the sphere. For
this, we check for which points, the gradients of the coordinate functions
and the position vector are linearly dependent. In other words, for which
points $p$ the following matrix has rank less than $3$:
$$ M = \left(
\begin{array}{cccccccc}
\frac{\partial \psi_1}{\partial x_1} & \frac{\partial \psi_1}{\partial y_1} & \cdots & \frac{\partial \psi_1}{\partial x_{n-1}} & \frac{\partial \psi_1}{\partial y_{n-1}} & \frac{\partial \psi_1}{\partial x_{n}}  & \frac{\partial \psi_1}{\partial y_{n}} \\
0  & 0  & \cdots & 0 & 0 &  \frac{\partial \psi_2}{\partial x_{n}} & \frac{\partial \psi_2}{\partial y_{n}} \\
x_1 & y_1 & \cdots & x_{n-1} & y_{n-1} & x_n & y_n\\
\end{array}
\right),$$
since $\frac{\partial \psi_2}{\partial x_j}=\frac{\partial \psi_2}{\partial y_j}=0$ for every $j=1, \cdots, n-1$.

To simplify notation we define the column vectors
$$ C(x_j) = \left(
\begin{array}{c}
\frac{\partial \psi_1}{\partial x_j} \\
 \frac{\partial \psi_2}{\partial x_{j}} \\
 x_j\\
\end{array}
\right), \ \ \
C(y_j) = \left(
\begin{array}{c}
\frac{\partial \psi_1}{\partial y_j} \\
 \frac{\partial \psi_2}{\partial y_{j}} \\
 y_j\\
\end{array}
\right).$$

Now consider the different types of minor of the matrix $M$. Note that
the interesting minors must contain at least one of the columns $C(x_n)$ or $C(y_n)$.
It is easy to see that the minors
$$|C(x_j ) C(y_j) C(x_n)| = |C(x_j) C(y_j) C(y_n)| = 0,$$
for $j = 1,\cdots,n-1$. So we do not need to consider them.
Now let us for $j = 1,\cdots,r$ the minors of the form:
\begin{align}
|C(x_j) C(x_n) C(y_n)| = \ \ \ \ \ \  \ \ \ \ \ \ \ \ \ \ \ \ \ \ \ \ \ \ \ \ \ \ \ \ \ \ \ \ \  \ \ \ \ \ \ \ \ \ \ \ \ \ \ \ \ \ &  \nonumber \\
 = \left|
\begin{array}{cccc}
2 \mu_ja_jx_j(x_j^2+y_j^2)^{a_j-1} & 3\mu_nx_n^2+\mu_ny_n^2 & 2\mu_nx_ny_n \\
0 &  2\mu_nx_ny_n & 3\mu_ny_n^2+\mu_nx_n^2 \\
x_j & x_n &  y_n\\
\end{array}
\right| & \nonumber \\
= x_j\mu_n(x_n^2+y_n^2)(3\mu_n(x_n^2+y_n^2)-2\mu_ja_j(x_j^2+y_j^2)^{a_j-1}x_n). &\end{align}

Analogously we have
\begin{align}
|C(y_j) C(x_n) C(y_n)| = \ \ \ \ \ \ \ \ \ \ \ \ \ \ \ \ \ \ \ \ \ \ \ \ \ \ \ \ \ \ \ \ \ \ \ \ \ \ \ \ \ \ \ \ \ \ \ \ & \nonumber \\
 = y_j\mu_n(x_n^2+y_n^2)(3\mu_n(x_n^2+y_n^2)-2\mu_ja_j(x_j^2+y_j^2)^{a_j-1}x_n). &
\end{align}

There are three cases when
$$|C(x_j) C(x_n) C(y_n)| = |C(y_j) C(x_n) C(y_n)| = 0.$$

\begin{itemize}
  \item [1.] If $x_j = y_j = 0$ for $j = 1,\cdots,n-1$, this gives points of the form
\begin{align}\label{eq25}
(0,0,\cdots,0,0,x_n,y_n).
\end{align}

  \item[2.] If $x_n = y_n = 0$, this gives points of the form
\begin{align}\label{eq26}
(x_1,y_1,\cdots,x_{n-1},y_{n-1},0,0).
\end{align}
  \item[3.] If $x_n^2 +y_n^2 \neq 0, x_j^2 +y_j^2 \neq 0$ for $j=1, \cdots, r$ and they satisfy the equation
\begin{align}\label{eq111}
3\mu_n(x_n^2+y_n^2)-2\mu_ja_j(x_j^2+y_j^2)^{a_j-1}x_n=0.
\end{align}
\end{itemize}

Note that in this case from (\ref{eq111}) if $x_n=0$ implies $y_n=0$, so in this case $x_n \neq 0$. Hence, equation (\ref{eq111}) is
equivalent to
\begin{align}\label{eq1234}
\dfrac{3\mu_n(x_n^2+y_n^2)}{x_n} = 2\mu_ja_j(x_j^2+y_j^2)^{a_j-1},
\end{align}
which can be written in the complex coordinates as
$$\dfrac{3\mu_n||z_n||^2}{Re(z_n)}= 2 \mu_ja_j||z_j||^{2a_j-2}.$$
Analogously, for $k=r+1, \cdots, n-1$ we have
\begin{align}\label{eq12345}
|C(x_k) C(x_n) C(y_n)| = \ \ \ \ \ \ \ \ \ \ \ \ \ \ \ \ \ \ \ \ \ \ \ \ \ \ \ \ \ \ \ \ \ \ \ \ \ \ \ \ \ \ \ \ \ \ \ \ & \nonumber \\
=  x_k\mu_n(x_n^2+y_n^2)(3\mu_n(x_n^2+y_n^2)+2\mu_ka_k(x_k^2+y_k^2)^{a_k-1}x_n). & \\
|C(y_k) C(x_n) C(y_n)| = \ \ \ \ \ \ \ \ \ \ \ \ \ \ \ \ \ \ \ \ \ \ \ \ \ \ \ \ \ \ \ \ \ \ \ \ \ \ \ \ \ \ \ \ \ \ \ \ &  \nonumber \\
 = y_k\mu_n(x_n^2+y_n^2)(3\mu_n(x_n^2+y_n^2)+2\mu_ka_k(x_k^2+y_k^2)^{a_k-1}x_n). &
\end{align}
These two minors are zero if $x_n^2+y_n^2 \neq 0, x_k^2 +y_k^2 \neq 0$ for $k=r+1, \cdots, n-1$ and they satisfy the equation
\begin{align}\label{eq123456}
 3\mu_n(x_n^2+y_n^2)-2\mu_ka_k(x_k^2+y_k^2)^{a_k-1}x_n=0.
\end{align}
Again, if $x_n=0$ then $y_n=0$, so in this case $x_n \neq 0$. Hence, equation (\ref{eq123456}) is equivalent to
\begin{align}\label{eq90}
\dfrac{3\mu_n(x_n^2+y_n^2)}{x_n}= 2\mu_ka_k(x_k^2+y_k^2)^{a_k-1},
\end{align}
which can be written in the complex coordinates as
$$\dfrac{3\mu_n||z_n||^2}{Re(z_n)}=-2\mu_ka_k||z_k||^{2a_k-2}.$$

\noindent \textbf{Claim}: Suppose $x_n^2 +y_n^2 \neq 0, x_j^2 +y_j^2 \neq 0$ for $j=1, \cdots, r$ and the satisfy equation (\ref{eq1234}).
Then $x_n$ and $y_n$ do not satisfy equation (\ref{eq90}) for any $x_k$ and $y_k$ such that $x_k^2 +y_k^2 \neq 0$ with $k=r+1, \cdots, n-1$.

\vspace{0.2cm}

\noindent\emph{Proof of Claim}: The right hand side of (\ref{eq1234}) is positive, then the left hand side must be also positive and this implies that $x_n>0$. Now suppose that $x_n$ and $y_n$ satisfy equation (\ref{eq90}) for some $x_k$ and $y_k$ such that $x_k^2 +y_k^2 \neq 0$ with $k=r+1, \cdots, n-1$. The right hand side of equation (\ref{eq90}) is negative, then the left hand side should also be negative and this implies that $x_n$ is negative. Then $x_n=0$, but both equations (\ref{eq1234}) and (\ref{eq90}) imply that $y_n =0$, which is a contradiction with the hypothesis that $x_n^2 +y_n^2 \neq 0$. \hspace{3.6cm} $\Box$

\vspace{0.1cm}
By the Claim we have that if $x_n$ and $y_n$ satisfy one of the equations (\ref{eq1234}) or (\ref{eq90}) it cannot satisfy the order, hence in principle we have points of the form
\begin{align}\label{eq33}
(x_1,y_1, \cdots, x_r,y_r,0,0,\cdots,0,0,x_n,y_n)
\end{align}
or
\begin{align}\label{eq34}
(0,0, \cdots,0,0,x_{r+1},y_{r+1},\cdots,x_{n-1},y_{n-1},x_n,y_n).
\end{align}

The rest of the minors that we need to analyze are of the forms
$$|C(x_j) C(x_k) C(x_n)|, \ \ \ \ |C(x_j) C(x_k) C(y_n)|, \ \ \ \ |C(y_j) C(y_k) C(x_n)|,$$
$$|C(y_j) C(y_k) C(y_n)|, \ \ \ \ |C(x_j) C(y_k) C(x_n)|, \ \ \ \ |C(x_j) C(y_k) C(y_n)|.$$

Is easy to see the following:
\begin{enumerate}
  \item[1.] In point of the form (\ref{eq25}) all the minors are zero and it is easy to see that points of this form cannot approach $V$.
  \item[2.] In point of the form (\ref{eq26}) all the minors are zero, and since $j = 1, \cdots, n-1$ are the critical indices of $\psi$, by Proposition \ref{proppontoscriticos} the points of the form (\ref{eq26}) are the critical points of $\psi$.
  \item[3.] If $x_n^2 +y_n^2 \neq 0$ we need to consider points of the forms (\ref{eq33}) and (\ref{eq34}).

  We have some cases
  \begin{enumerate}
    \item[(a)] If one of the above minors contains the columns $C(u_j)$ and $C(u_k)$ with $u=x$ or $u=y$ with $j=1, \cdots,n-1$, then both $u_j$ and $u_k$ are factors, and since one of them has to be zero, the minor is zero.
    \item[(b)] If one of the above minors contains the columns $C(u_j)$ and $C(u_k)$ with $u=x$ or $u=y$ with $j,k=1, \cdots, r$(or $j,k=r+1, \cdots, n-1$) they need to satisfy the equation
     \begin{align}\label{eq35}
   (2\mu_ja_j(x_j^2 +y_j^2)^{a_j-1} - 2\mu_k a_k(x_k^2 +y_k^2)^{a_k-1})=0.
     \end{align}
     But since $x_n$ and $y_n$ satisfy equations (\ref{eq1234}) for $j$ and $k$ we have
     $$\dfrac{3 \mu_n(x_n^2 +y_n^2)}{x_n}=2\mu_ja_j(x_j^2 +y_j^2)^{a_j-1},$$
     and
     $$\dfrac{3 \mu_n(x_n^2 +y_n^2)}{x_n}=2\mu_ka_k(x_k^2 +y_k^2)^{a_k-1},$$
     (for $j,k=r+1, \cdots,n-1$ there is a minus sign in the right hand side of both equalities) thus, both terms in equation (\ref{eq35}) are equal and its difference is zero. Therefore all the minors are zero in points of the form (\ref{eq33}) and (\ref{eq34}).
     Now we need to check that this type of points cannot approach $V$. Recall that for points of type (\ref{eq33}) we have $x_n > 0$ and for points of type (\ref{eq34}) we hhave $x_n < 0$. We check the first case, the second one is analogous. Suppose there is a sequence of points
     $$p_l=(x_{1,l},y_{1,l}, \cdots, x_{r,l},y_{r,l},0,0,\cdots, 0,0, x_{n,l},y_{n,l}),$$
     of type (\ref{eq33}) converging to a point of the form
     $$p=(x_1,y_1, \cdots, x_r,y_r,0,0,\cdots, 0,0,x_n,0),$$
  that is, with the last coordinate zero. For $p$ to be in $V$ in has to satisfy
  $$\sum_{j=1}^{r}\mu_j(x_j^2 +y_j^2)^{a_j} + \mu_n x_n^3 =0.$$

  Since $\sum_{j=1}^{r}\mu_j(x_j^2 +y_j^2)^{a_j} > 0$ we need $x_n^3 <0$, i.e., $x_n < 0$, but since the points $p_l$ are of the form (\ref{eq33}) we have $x_{n,l}>0$ for all $l=0,1,\cdots$, thus $x_n \geq 0$ and $p$ cannot be in $V$.
  \end{enumerate}
\end{enumerate}

By Proposition \ref{teopropdetransversalidade}, $\psi$ described by equations (\ref{eq21}) and (\ref{eq12}) has the transversality property.

We know that $\psi$ is $d$-regular (Proposition \ref{propfibracaodafamiliad-regular}), so by the Theorem \ref{teodafibracaoCMSS} we prove that $\psi$ has Milnor Fibration.
\end{prova}

\begin{remark}

Note that Lemma \ref{lemacasoespecial} is valid even if all $\mu_j$ are positive or all are negative, so these particular cases are also covered by Theorem \ref{teofibracaoitens}.
\end{remark}

\begin{theorem}\label{teoremaespecial}
Consider mixed polynomials given by
\begin{equation}\label{eq36}
\psi(z)=\sum_{j=1}^{n-1}\lambda_j z_j^{a_j}\overline{z_j}^{a_j} + \lambda_n z_n^2 \overline{z_n}
\end{equation}
or
\begin{equation}\label{eq37}
\psi(z)=\sum_{j=1}^{n-1}\lambda_j z_j^{a_j}\overline{z_j}^{a_j} + \lambda_n z_n \overline{z_n}^2,
\end{equation}
with $a_j, b_j \in \mathbb{Z}, a_j,b_j >0$ and all the $\lambda_j$ for $j=1, \cdots, n$ colinear. Then the map
$$\phi=\dfrac{\psi}{||\psi||}:(\mathbb{S}_{\epsilon}^{2n-1} \setminus W) \to \mathbb{S}^1 \setminus \mathcal{A},$$
is a locally trivial fibration, where $\Delta_\psi^\eta= \Delta_\psi \cap \mathbb{B}_\eta^2, W=\psi^{-1}(\Delta_\psi^\eta), \mathcal{A}_\eta = \Delta_\psi^\eta \cap \mathbb{S}^1_\eta$ and $\mathcal{A}=\pi(\mathcal{A}_\eta)$ and $\pi: \mathbb{C}^* \to \mathbb{S}^1$ is the radial projection.
\end{theorem}
\vspace{-0.2cm}
\begin{prova}
We write the proof for (\ref{eq36}), the other is analogous. Let $\theta=\arg \lambda_1$. Since all $\lambda_j$ with $j=1, \cdots, n$ are colinear we have that $\lambda_j=\mu_j e^{i \theta}$ with $\mu_j \in \mathbb{R}$ for $j=1, \cdots, n$. Hence we have
$$\psi(z)= \left(\sum_{j=1}^{r}\mu_jz_j^{a_j}\overline{z_j}^{a_j} + \mu_nz_n^2 \overline{z_n}\right)e^{i\theta}.$$

Let $r_\theta: \mathbb{C} \to \mathbb{C}$ be the analytic map given by multiplication by $e^{-i\theta}$. Let $\tilde{\psi}$ be the composition $\tilde{\psi}= r_\theta \circ \psi$, thus we have
$$\tilde{\psi}(z)= \sum_{j=1}^{r}\mu_j z_j^{a_j}\overline{z_j}^{a_j} + \mu_n z_n^2 \overline{z_n}.$$

Up to renumbering of the variables (so that $\mu_j > 0$ for $j = 1,\cdots,r$ and $\mu_j < 0$ for $j = r +1,\cdots , n-1)$ we have that $\tilde{\psi}$ is of the form (\ref{eq21}) or (\ref{eq12}) in Lemma \ref{lemacasoespecial}. Therefore $\tilde{\psi}$ has a fibration on the sphere and since $r_\theta$ is a diffeomorphism also $\psi$ has a fibration on the sphere.

\end{prova}

\section*{Acknowledgments}

The authors are grateful to Jos\'e Seade for their careful reading of the first draft of this paper and for their suggestions, and Maria Ruas, Jos\'e Luis Cisneros-Molina, Jawad Snoussi and Aur\'elio Menegon Neto for many discutions and suggestions which greatly improve the content and presentation of this work.

This work was developed during the doctorate of the second author at ICMC-USP (Brazil) and concluded during the Thematic Program on Singularities at IMPA (Brazil), February 2020.

The first author was supported by CNPq (Conselho Nacional de \linebreak Desenvolvimento Cient\'ifico e Tecnol\'ogico, Brazil) Grant $303046/2016-3$ and FAPESP (Funda\c c\~ao de Amparo \`a Pesquisa do Estado de S\~ao Paulo, Brazil) Grant $2017/09620-2$.

The secound author was supported by FAPESP (Funda\c c\~ao de Amparo \`a Pesquisa do Estado de S\~ao Paulo, Brazil) Grant $2013/22718-0$, $2015/19721-5$.

We thank the referee for his/her carefully reading, remarks and suggestions that provided generalizations and elegance for the main results of this work.

\end{document}